\font\titlefont= cmcsc10 at 12pt
\def\A{{\bf A}}
\def\F{{\bf F}}
\def\P{{\bf P}}
\def\Z{{\bf Z}}

\font\sixteenbf=cmbx10 scaled\magstep2
\def\vandaag{\number\day\space\ifcase\month\or
 januari\or februari\or  maart\or  april\or mei\or juni\or  juli\or
 augustus\or  september\or  oktober\or november\or  december\or\fi,
\number\year}
\magnification\magstep1
%
%
\catcode`\@=11
\font\tenmsa=msam10
\font\sevenmsa=msam7
\font\fivemsa=msam5
\font\tenmsb=msbm10
\font\sevenmsb=msbm7
\font\fivemsb=msbm5
\newfam\msafam
\newfam\msbfam
\textfont\msafam=\tenmsa  \scriptfont\msafam=\sevenmsa
  \scriptscriptfont\msafam=\fivemsa
\textfont\msbfam=\tenmsb  \scriptfont\msbfam=\sevenmsb
  \scriptscriptfont\msbfam=\fivemsb
\def\hexnumber@#1{\ifcase#1 0\or1\or2\or3\or4\or5\or6\or7\or8\or9\or
        A\or B\or C\or D\or E\or F\fi }
\edef\msa@{\hexnumber@\msafam}
\edef\msb@{\hexnumber@\msbfam}
\mathchardef\square="0\msa@03
\mathchardef\subsetneq="3\msb@28
\mathchardef\ltimes="2\msb@6E
\mathchardef\rtimes="2\msb@6F
\def\Bbb{\ifmmode\let\next\Bbb@\else
        \def\next{\errmessage{Use \string\Bbb\space only in math mode}}\fi\next}
\def\Bbb@#1{{\Bbb@@{#1}}}
\def\Bbb@@#1{\fam\msbfam#1}
\catcode`\@=12
%
%

%
%
\vglue 2.0cm
\centerline{\sixteenbf The Coset Weight Distributions of Certain BCH Codes} 
\bigskip
\centerline{\sixteenbf  and a Family of Curves}
\bigskip
\vskip 2pc
\centerline{\titlefont G.\ van der Geer \& M.\  van der Vlugt}
\bigskip\bigskip
\centerline{\bf  Introduction}
\bigskip
\noindent
Many problems in coding theory are related to the problem of
determining the distribution of the number of rational points
in a family of algebraic curves defined over a finite field.
Usually, these problems are very hard and a complete answer is
often out of reach.

In the present paper we consider the problem of the weight distributions
of the cosets of certain BCH codes. This problem turns out to be
equivalent to the determination of the distribution of 
the number of points in a family of curves 
with a large symmetry group. The symmetry allows us to 
analyze  closely the nature of these curves and in this way we are able to
extend considerably our control over the coset
weight distribution compared with earlier results.

For a binary linear code $C$ of length $n$ the weight distributions of the
cosets of $C$ in $\F_2^n$ are important invariants of the code.
They  determine for example the probability 
of a decoding error when using $C$. However, the coset weight 
distribution problem is solved for very few types of codes.

In [C-Z] Charpin and Zinoviev study the weight distributions of the cosets
of the binary $3$-error-correcting BCH code of length $n=2^m-1$ with $m$ odd.
We denote this code by $BCH(3)$.

Let $\F_q$ be a finite field of cardinality $q=2^m$ and let 
$\alpha$ be a generator of the multiplicative group $\F_q^*$. The matrix
$$
H= \left(\matrix{ 1 & \alpha & \alpha^2 & \ldots & \alpha^{n-1}\cr
1 & \alpha^3 & \alpha^6 & \ldots & \alpha^{3(n-1)}\cr
1 & \alpha^5 & \alpha^{10} & \ldots & \alpha^{5(n-1)}\cr
}\right)
$$
is a parity check matrix defined over $\F_q$ of $BCH(3)$. This means that
$$
BCH(3)=\{  c=(c_0,\ldots,c_{n-1})\in \F_2^n : H c^t=0 \}.
$$

It was shown in [C-Z] that the coset weight distribution problem
for $BCH(3)$ comes down to the same problem for the extended code
$\widehat{BCH(3)}$ with parity check matrix
$$
\hat{H}= \left(\matrix{ 1& 1 & 1& \ldots & 1&1\cr
1 & \alpha & \alpha^2 & \ldots & \alpha^{n-1}&0\cr
1 & \alpha^3 & \alpha^6 & \ldots & \alpha^{3(n-1)}&0\cr
1 & \alpha^5 & \alpha^{10} & \ldots & \alpha^{5(n-1)}&0\cr
}\right)
$$
A coset $\hat{D}$ of $\widehat{BCH(3)}$ in $\F_2^{n+1}$   
is characterized by
the syndrome $s(\hat{D})=\hat{H}x^t \in \F_q^4$, 
where $x$ is a representative of $\hat{D}$.
The weight of $\hat{D}$ is the minimum weight of the vectors in $\hat{D}$. Here
the weight of a vector is the number of its non-zero entries.

Charpin and Zinoviev then show that the weight distribution problem 
for the cosets of $\widehat{BCH(3)}$ of length $2^m$ with $m$ odd can be
solved as soon as the weight distributions of the cosets $\hat{D_4}$
of weight $4$ with syndrome $s(\hat{D_4})=(0,1,A,B)$ are determined.

From [C-Z] we recall:
{\sl The weight distribution of a coset $\hat{D_4}$
is determined by the number $N(A,B)$ of vectors of weight $4$ in $\hat{D_4}$.}

Via the matrix $\hat{H}$ this leads to the  system of equations in four 
variables in $\F_{q=2^m}$:
$$
\eqalign{
x_1+x_2+x_3+x_4 &= 1,\cr
x_1^3+x_2^3+x_3^3+x_4^3 &= A,\cr
x_1^5+x_2^5+x_3^5+x_4^5 &= B, \cr
} \eqno(1)
$$
and $N(A,B)$ is the number of $S_4$-orbits of solutions of (1)
with distinct $x_i \in \F_q$. In particular the number of values
of $N(A,B)>0$ equals the number of different coset weight distributions
of cosets of type $\hat{D_4}$. Note that since the set of solutions of (1)
is invariant under translation over $(1,1,1,1)$ the quantity
$N(A,B)$ is even.

In this paper we shall show that by analyzing carefully the curves
defined by (1) we can determine good upper and lower bounds
for the pivotal quantity $N(A,B)$. The bounds are obtained 
by dissecting the Jacobian variety of the curves in our family 
in isogeny factors of dimension $1$ and $2$. This yields restrictions
on the traces of Frobenius. The splitting of the Jacobian is a
corollary from a very effective  description of the curves defined by (1) as
fibre products over $\P^1$ of three elliptic curves.

We show that for odd $m$ the $N(A,B)$ lie in an explicit interval of length
$\sim 1.57 \sqrt{q}$, cf.\ 
[C-Z], where the interval is $\sim q/4$. 
Moreover, we argue that on statistical grounds
one may expect that almost all $N(A,B)$ lie in an explicit interval
of length $\sim 0.9 \sqrt{q}$. We then give
numerical results that confirm strongly these heuristics and extend
the table of $BCH(3)$ codes with known coset weight distribution.

\bigskip
\noindent
\centerline{\bf \S 1. A family of curves}
\bigskip
\noindent
We consider the algebraic curve $C^{\prime}=C^{\prime}_{A,B}$ in $\P^4$
given by the equations
$$
s_1 = x_0,\quad
s_3 = Ax_0^3,\quad
s_5 = Bx_0^5, \eqno(2)
$$
where  $s_j$ is the $j$-th power sum $\sum_{i=1}^4 x_i^j$ in the variables
$x_1,\ldots,x_4$. Let $\sigma_j$ denote the $j$-th elementary
symmetric function in $x_1,\ldots,x_4$. 
If we apply Newton's formulas for power sums we find
$$
\eqalign{ 
s_1+x_0 &= \sigma_1+x_0=0,\cr
s_3+Ax_0^3& = (A+1)x_0^3 +\sigma_2x_0+\sigma_3=0,\cr
s_5+Bx_0^5 &= x_0\, ((B+A)x_0^4+(A+1)\sigma_2x_0^2+\sigma_4)=0.\cr
}$$
This implies that 
the curve $C^{\prime}$ consists of the three lines in the hyperplane
$x_0=0$ given by
$$
x_i+x_j=x_k+x_l=0, \quad \hbox{\rm with} \quad \{i,j,k,l\}=\{1,2,3,4\},
\eqno(3)
$$
and a curve $C=C_{A,B}$ given by
$$
\eqalign{
\sigma_1 &= x_0,\cr
\sigma_3 &= (A+1)x_0^3+\sigma_2x_0,\cr
\sigma_4 &= (B+A)x_0^4+(A+1)\sigma_2 x_0^2. \cr}\eqno(4)
$$

The symmetric group $S_4$ operates on $C^{\prime}$ and on $C$ by
permuting the coordinates $x_1,\ldots,x_4$.
Moreover, there is an involution $\tau$ acting on $C$ via
$$
(x_0:x_1:\ldots:x_4) \mapsto (x_0:x_1+x_0:\ldots :x_4+x_0).
$$
This involution commutes with the elements of $S_4$ and this gives
rise to a group $G$ of $48$ automorphisms of $C$.

We introduce the invariant
$$
\lambda:= B+A^2+A+1 \quad (\in \F_q).
$$
\proclaim (1.1) Lemma. 
\item{i)}  If $\lambda\neq 0$ then $C$ has six ordinary double points,
namely the points of the $S_4$-orbit of $(0:1:1:0:0)$ and no other
singularities.
\item{ii)} If $\lambda=0$ the curve $C$ consists of $12$ lines.
\par
\noindent
{\sl Proof.} 
The Jacobian matrix of (2) is
$$
\left( \matrix{1&1&1&1&1\cr
Ax_0^2&x_1^2&x_2^2&x_3^2&x_4^2\cr
Bx_0^4&x_1^4&x_2^4&x_3^4& x_4^4\cr} \right) .
$$
If the rank of this matrix is $\leq 2$ for a point with coordinates 
$(x_0:\ldots : x_4)$ then there exist $\alpha, \beta, \gamma$
with $\alpha$, $\beta$, $\gamma$ not all zero 
 such that $\alpha + \beta x_i^2+\gamma x_i^4=0$
for $i=1,\ldots,4$. Hence the coordinates $x_i$ with $i=1,\ldots,4$ 
of a singular point of $C$
can assume at most $2$ different values and 
taking into account the equation $s_1=x_0$ 
it follows that a singular point of $C$ is 
in the $S_4$-orbit of a point of the form 
$(a:1:1:1:a+1)$ or of the form $(0:1:1:a:a)$ for some value of $a$.
In the latter case we get from (4) that $a=0$ and we find $6$ 
singular points in the orbit of $(0:1:1:0:0)$. 
In the former case it follows from (4) that $a$ satisfies
$$
(A+1)a^3+a^2+a=0,\quad \hbox{\rm and} \quad
(B+A)a^4+(A+1)a^3+a+1=0. \eqno(5)
$$
Hence $a\neq 0$ and (5) is equivalent with
$$
\eqalign{
(A+1)a^2+a+1 &=0, \cr
(B+A)a^2+(A+1)a+(A+1) &=0.\cr
} \eqno(6)
$$
The resultant of (6) equals $(B+A^2+A+1)^2$, hence (6) has a solution
if and only if $\lambda=B+A^2+A+1$  vanishes. 
In that case the Jacobian matrix has rank $2$ for the solutions of $(6)$.

So if $\lambda \neq 0$ the curve $C$ has six singular points,
namely the $S_4$-orbit of $(0:1:1:0:0)$. For the local structure 
near $(0:1:1:0:0)$ we eliminate $x_0$ from (2) and find that the curve
$C^{\prime}$ in $\P^3$ is given by 
$$
s_3= As_1^3, \quad s_5=Bs_1^5.
$$
Upon taking affine coordinates
$\xi_1=(x_1+x_2)/x_1$, $\xi_2=x_3/x_1$, 
$\xi_3=x_4/x_1$ we find the equations
$$ 
\eqalign{
\xi_1+\xi_1^2+\xi_1^3+\xi_2^3+\xi_3^3& = A(\xi_1+\xi_2+\xi_3)^3, \cr
 \xi_1+\xi_1^4+\xi_1^5+\xi_2^5+\xi_3^5& = B(\xi_1+\xi_2+\xi_3)^5. \cr
}
$$
This shows that $\xi_1$ lies in $m^3$, with $m$ the maximal ideal of $(0,0,0)$
in $\A^3$ and defines the tangent plane at $(0,0,0)$ 
to the cubic surface $S$ given by
the cubic equation. Moreover, 
this is also the lowest order term of the quintic 
equation. Therefore, locally near the origin $C^{\prime}$ is given by
$$
\xi_1=0, \quad (\xi_2+\xi_3)(\xi_2\xi_3 + (A+1)(\xi_2+\xi_3)^2)=0.\eqno(7)
$$
which shows that $C^{\prime}$ has a triple point and 
$C$ has a node at this point.

If $\lambda=0$ and $a$ 
satisfies
$
(A+1)a^2+a+1=0
$
then  $a$ is a solution of (6) and the $S_4$-orbit of points of the 
form $(a:x:x:1:a+1)$ with arbitrary $x$ is on $C$. So the equations
$$
x_i+x_j=0, \quad (a+1)x_k+x_l=0 \quad \hbox{\rm with} \quad \{i,j,k,l\}=\{1,2,3,4\}
$$
define a line on $C$ and this gives $12$ lines on $C$. Since
$C$ has degree $12$ the curve $C$ decomposes as the union of
$12$ lines. This proves ii).
\bigskip
\noindent
{\bf Remark}. 
It follows from the preceding proof that for $\lambda \neq 0$ points on $C$
for which $x_1,\ldots,x_4$ are not all distinct
lie on one of the lines (3).
\medskip

\proclaim  (1.2) Proposition. If $\lambda\neq 0$ then $C$ is irreducible. \par
\noindent
{\sl Proof.} Suppose that $C=\sum_{i=1}^m C_i$ is a sum of 
irreducible components $C_i$
with $m\geq 2$. Since $C$ is connected at least one of the singular points
is an intersection point of two distinct components $C_i$. 
By the $S_4$-symmetry
then each of the six singular points is an intersection point 
of two different components. 
This implies that the components $C_i$ are non-singular.
Since the permutation $(34)$ interchanges the two branches of $C$
in $(0:1:1:0:0)$ (cf.\ (7)) the group $S_4$ acts transitively on the branches
through a singular point, so $S_4$ acts transitively on the set of
components.

Let $S$ be the smooth cubic surface in $\P^4$ 
given by the equations $s_1=x_0, \, s_3=Ax_0^3$.
On $S$ the curve $C$ is linearly equivalent to $4H$ with $H$ the hyperplane
section of $S$. Now the intersection number $HC_i$ equals the intersection
number with the hyperplane $x_0$, i.e.\ the intersection number of $C_i$
with the three lines (3), and since the intersection is transversal 
$HC_i$ equals the number of singular points 
of $C$ on $C_i$. Put $r= 12/m$. Then by the symmetry 
we have $HC_i=r$. On the other hand,
the adjunction formula
$$
C_i^2+K_SC_i=C_i^2-HC_i= C_i^2-r=2g(C_i)-2,
$$
where $K_S$ is the canonical divisor of $S$,  and  the identity
$$
4r= 4HC_i= CC_i=C_i^2+\sum_{j\neq i} C_iC_j=C_i^2+r
$$
imply $C_i^2=3r$ and $g(C_i)=r+1$. In particular, $C_i$ cannot be
contained in a hyperplane and spans $\P^3$.
Clifford's theorem applied to the hyperplane section $H|C_i$ of $C_i$
says that $h^0(H|C_i)\leq r/2+1$, hence $r\geq 6$. Then $m=2$ and we have
two components. Again, by Clifford, these curves must be hyperelliptic and 
the linear system $H|C_i$ is $3g_2^1$. But since $3g_2^1$ is contained 
in the canonical system $|K_{C_i}|$ this factors through the hyperelliptic 
involution, which contradicts the fact that $C_i$ is embedded in $\P^3$ as
a non-rational curve. This proves that $C$ is irreducible.

\proclaim (1.3) 
Corollary. If $\lambda \neq 0$ the normalization 
 $\tilde{C}$ of $C$ is an irreducible 
smooth curve of genus $13$. \par
\noindent
{\sl Proof.} On the cubic surface $S$ we have
$ (C+K_S)C=(4-1)HC=36$. 
This implies that for $\tilde{C}$ we have 
$2g(\tilde{C})-2=36-12=24$. $\square$
\bigskip
\bigskip
\noindent
\centerline{\bf \S 2. Dissecting the Jacobian}
\bigskip
\noindent
We now analyze the curve $C=C_{A,B}$ for 
$\lambda \neq 0$ in more detail in order
to decompose its Jacobian.
\medskip
Let $H\subset {\rm Aut}(C)$ be the subgroup generated by the two 
permutations $(12)$ and $(34)$ and the involution $\tau$.
Then $H$ is abelian of order $8$ and isomorphic to $(\Z/2\Z)^3$.
Consider the following diagram of degree $2$ coverings of curves
$$
\matrix{
&&& C &&&\cr
&&\swarrow & \downarrow & \searrow&\cr
& C/(12) && C/(12)(34) && C/(34)&\cr
&\downarrow & \searrow & \downarrow & \swarrow & \downarrow &\cr
& C/\langle (12),\tau\rangle && C/\langle (12),(34)\rangle && 
C/\langle (34),\tau \rangle &\cr
&&\searrow & \downarrow & \swarrow && \cr
&&& C/H &&&\cr
}
$$

Let $u = (x_3+x_4)/(x_1+x_2)$. This is a $H$-invariant rational function
on $C$, hence defines a rational function on $C/H$.

\proclaim (2.1) Proposition. \item{i)} The function $u$ gives an isomorphism
$C/H \cong \P^1$.
\item{ii)} The curve $C/\langle (12),(34)\rangle$ is a curve of genus $1$
given by
$$
y^2+y= \lambda u + \lambda /u +(A+1).
$$
\item{iii)} The curve
 $C/\langle (12), \tau)\rangle$ is a curve of genus $2$ given by
$$
y^2+y=\lambda /u^3+\lambda u.
$$
\item{iv)}  The curve
$C/\langle (34), \tau\rangle$ is a curve  of genus $2$ given by
$$
y^2+y=\lambda u^3+\lambda /u.
$$
\par
\noindent
{\sl Proof.} 
The divisor of $u$ on $C$ is of the form $H_{34}C -H_{12}C$,
where $H_{ij}$ is the hyperplane given by $x_i+x_j=0$. 
Since both these hyperplanes contain the line $x_1+x_2=x_3+x_4=0$ 
which intersects $C$ in a divisor of degree $4$
it follows that the divisor of $u$
can be written as a difference of two divisors of degree $12-4=8$. Moreover,
these divisors are invariant under the action of $H$. This implies that on
$C/H$ the function $u$ defines a non-constant function with a single
zero and a single pole. Therefore $u$ 
defines an isomorphism $C/H \cong \P^1$. This proves i).

We now prove ii). Working with the affine equations (set $x_0=1$)
$$
\sigma_1 =1, \quad
\sigma_3 = A+1+\sigma_2, \quad
\sigma_4=B+A+ (A+1)\sigma_2, 
$$
we can write $u= (x_3+x_4)/(x_1+x_2)= 1+ 1/(x_1+x_2)$,
i.e. 
$$
x_1+x_2=1/(u +1)\quad
\hbox{\rm  and }\quad x_3+x_4= u/ (u +1).
$$
We put $v:= x_1x_2$ and  $w:= x_3x_4$. These functions are invariant
under $(12)$ and $(34)$, but not under $\tau$.
Using 
$$
\eqalign{
\sigma_2 &= x_1x_2 + x_3x_4+ (x_1+x_2)(x_3+x_4)
         = v+w+u/(u+1)^2,\cr
\sigma_3 &= x_1x_2(x_3+x_4) + (x_1+x_2)x_3x_4
= (uv +w)/(u +1),\cr
}
$$
the equation $\sigma_3=A+1+\sigma_2$ implies  
$$
A(u+1)^2 + (u+1)(v +uw) + u^2+u +1=0, \eqno(8)
$$
while the equation $\sigma_4= B+A+ (A+1)\sigma_2$ yields
$$
B+A + (A+1)(v+w+u/(u+1)^2) +vw=0. \eqno(9)
$$
Elimination of $w$ from (8) and (9) yields the equation
$$
(u +1)^2 v^2+u (u +1) v= 
\lambda u^3 +\lambda u  +(A+1)u + (A+1)^2u^2 + (A+1)^2.
$$
Dividing by $u^2$ and and replacing $(u+1)v/u$ by $\eta$ (i.e.\
$\eta= x_1x_2/(x_3+x_4)$) gives
$$
\eta^2+\eta= \lambda u + \lambda/u +
(A+1)/u + (A+1)^2/u^2+ (A+1)^2 \eqno(10)
$$
and this is Artin-Schreier equivalent to
$$
y^2+y= \lambda u + \lambda/u + A+1.
$$ 
Since $\eta$ is invariant under $(12)$ and $(34)$, 
but not under $\tau$ the equation (10) describes 
the degree $2$ cover $C/\langle (12),(34)\rangle$ of $C/H$.

For iii) we remark that the function 
field extension of $C/\langle (12), \tau \rangle$
over $C/H$ is generated by the function $z= x_3 + x_1x_2/(x_3+x_4)$.
Then $z+z^{(34)}= x_3+x_4= u/(u+1)$. Moreover,
$$ 
\eqalign{
z \cdot z^{(34)} &= x_3x_4 + x_1x_2 + (x_1x_2)^2/(x_3+x_4)^2\cr
& = w + v+ \eta^2\cr
&= A(u +1)/u + 1 + 1/u (u +1) + \eta +\eta^2,\cr
}
$$
where we used 
$w= A(u+1)/u + 1+ 1/u (u +1) + v/u$
obtained from (8) and $v+ v/u=\eta$.
By (10) this implies that $z$ satisfies the equation
$$
z^2+ { u\over u+1} z = {\lambda(u^4+u^3)+ (A^2+A)u^3+
\lambda(u^2+u)+ A u^2+ A^2+1\over u^2 (u +1)} \, .
$$
Dividing by $(u/(u+1))^2$ and 
replacing $(u+1)z/u$ by $\zeta$ gives the equation
$$
\zeta^2+\zeta= \lambda u + \lambda/u^3 + (A^2+A) + 1/u + 
1/u^4+ A/u^2+A^2/u^4.
$$
Via $\zeta \mapsto \zeta + A  +1/u+ (A+1)/u^2$
we get the Artin-Schreier equivalent form 
$$
\zeta^2+\zeta= \lambda u + \lambda /u^3.
$$
Part iv) is now obtained by applying the permutation $(13)(24)$. This 
changes $u$ into $u^{-1}$  and proves the result.
\smallskip
\noindent
\proclaim (2.2) Theorem. The curve $C$ is the fibre product over $\P^1$ with
affine coordinate $x$ of the three hyperelliptic curves     
$$
\eqalign{
y^2+y&= \lambda x^3+ \lambda/x,\cr
y^2+y&= \lambda/x^3+ \lambda x, \cr
y^2+y&= \lambda x + \lambda /x + A+1.\cr
}
$$
\par
\noindent
{\sl Proof.} This follows directly from the diagram and the preceding
proposition.

\bigskip
Note that equivalently, $C$ is the fibre product of 
the three curves $C_{f_i}$ of genus $1$ 
with affine equation $y^2+y=f_i$, where $f_i$
for $i=1,2,3$ is given by
$$
\eqalign{
f_1&= \lambda x^3+ \lambda x + A+1,\cr
f_2&= \lambda /x^3+ \lambda /x + A+1, \cr
f_3&= \lambda x + \lambda /x + A+1.\cr
} \eqno(11)
$$
This description allows us to dissect  the Jacobian of $C$.

\proclaim (2.3) Theorem.
The Jacobian of $C_{A,B}$ decomposes up to isogeny over $\F_q$
as a product of five supersingular elliptic curves, 
two ordinary elliptic curves and three $2$-dimensional factors of $2$-rank $1$. \par
\noindent
{\sl Proof.} From the description of $C=C_{A,B}$ as a fibre product we see
that ${\rm Jac}(C)$ decomposes
as a product of seven factors: three elliptic curves ${\rm Jac}(C_{f_i})$, 
two $2$-dimensional factors ${\rm Jac}(C_{f_1+f_3})$, ${\rm Jac}(C_{f_2+f_3})$,
 and two $3$-dimensional factors ${\rm Jac}(C_{f_1+f_2})$ and 
${\rm Jac}(C_{f_1+f_2+f_3})$. 
The $2$-rank of ${\rm Jac}(C_{f_i})$ is $0$ for $i=1,2$ and $1$ for $i=3$. 
The $2$-ranks of ${\rm Jac}(C_{f_1+f_3})$
and ${\rm Jac}(C_{f_2+f_3})$ are $1$ since these hyperelliptic curves
have two Weierstrass points.

The curve $C_{f_1+f_2+f_3}$ is a curve of genus $3$ 
defined by $y^2+y=\lambda(x^3+1/x^3)+A+1$
with automorphisms
$$ \rho: (x,y)\mapsto (1/x,y), \quad \sigma: (x,y)\mapsto (x, y+1).
$$
The quotient of $C_{f_1+f_2+f_3}$
under $\rho$ is the supersingular elliptic curve given by $y^2+y=
\lambda(z^3+z)+A+1$ with $z=x+1/x$. 
Moreover, the curve $C_{f_1+f_2+f_3}$
admits a non-constant map to the ordinary elliptic curve 
$y^2+y=\lambda(w+1/w)+A+1$ via $w=x^3$. So by Poincar\'e's complete
irreducibility theorem the Jacobian 
${\rm Jac}(C_{f_1+f_2+f_3})$ splits up to
isogeny into a product of three elliptic curves and has $2$-rank $1$
since it has $2$ ramification points.

Similarly, the quotient of $C_{f_1+f_2}$  by the automorphism $\rho$
is the supersingular elliptic curve $y^2+y=\lambda z^3$, while the quotient 
under $\rho\sigma$ is a curve of genus $2$ of $2$-rank $1$ 
defined by the equation
$y^2+y=\lambda z^3+1/z$. Collecting these results we obtain 
the theorem.
\medskip
For a smooth absolutely irreducible complete 
curve $X$ defined over a field $\F_q$ 
we shall denote the trace of Frobenius by $t(X)$, i.e.\
$t(X)=q+1-\#X(\F_q)$.

\proclaim (2.4) Corollary. For $q=2^m$ with $m$ odd the trace of
Frobenius of $C_{A,B}$ equals $2t(C_{f_1})+2t(C_{f_3})+2t(C_{f_1+f_3})+
t(C_{g_{\lambda}})$,
where $C_{g_{\lambda}}$ is the curve given by $y^2+y=g_{\lambda}$ 
with $g_{\lambda}=\lambda x^3+ 1/x$.
\par
\noindent
{\sl Proof.} The curves $C_{f_1}$ and $C_{f_2}$ are isomorphic via 
$x \mapsto 1/x$, so have the same trace of Frobenius. Since for
$q=2^m$ with $m$ odd the map $x \mapsto x^3$ is a bijection on $\F_q$
the curve $C_{f_1+f_2+f_3}$ given by $y^2+y=\lambda(x^3+1/x^3) + A+1$
and the ordinary factor of its Jacobian given by $y^2+y=
\lambda(w+1/w)+A+1$ have the same trace of Frobenius, and this is
$t(C_{f_3})$.  Moreover, since $C_{f_1+f_3}$ and $C_{f_2+f_3}$ are
isomorphic, we have $t(C_{f_1+f_3})=t(C_{f_2+f_3})$.
Similarly, the supersingular 
component of ${\rm Jac}(C_{f_1+f_2})$ given by
$y^2+y=\lambda z^3$ has the same trace of Frobenius as the rational
curve $y^2+y=\lambda z$, i.e.\ $0$. Therefore, the trace $t(C_{f_1+f_2})$
equals the trace of the genus $2$ quotient $C_{f_1+f_2}/\rho \sigma$,
and this is the curve $y^2+y=g_{\lambda}$.
\bigskip

We can interpret and augment the results obtained using the involution
$\tau$.  The involution $\tau$ acts without fixed points on $C$, hence by
 the Hurwitz-Zeuthen formula the genus of the quotient curve $C/\tau$
is $7$. The Jacobian ${\rm Jac}(C)$ decomposes up to an isogeny 
$$
{\rm Jac}(C) \sim {\rm Jac}(C/\tau) \times P,
$$
where $P$ is the Prym variety of $C \to C/\tau$, i.e.\ the identity
component of the norm map ${\rm Nm}: {\rm Jac}(C) \to {\rm Jac}(C/\tau)$.
Since the curves $C/\langle (12), \tau\rangle =C_{f_2+f_3}$, 
$C/\langle (34), \tau\rangle =C_{f_1+f_3}$ 
and $C/\langle (12)(34), \tau\rangle= C_{f_1+f_2}$ are quotients of $C/\tau$
and the fibre product
$C_{f_1+f_3}\times_{\P^1}C_{f_2+f_3}$ has genus $7$ it follows
readily that 
$$
C/\tau \cong C_{f_1+f_3}\times_{\P^1}C_{f_2+f_3}. 
$$
Note that the substitution $x \mapsto x/\lambda$
yields an isomorphism $C_{g_{\lambda^4}} \cong
C_{f_1+f_3}$.

\proclaim
(2.5) Proposition.
Up to isogeny over $\F_{q=2^m}$ we have the splitting
$$
{\rm Jac}(C/\tau) \sim {\rm Jac}(C_{g_{\lambda^4}})^2 
\times {\rm Jac}(C_{g_{\lambda}}) \times E,
$$
where $C_{g_{\lambda}}$ is as in (2.4) and
$E$ is the elliptic curve $y^2+y=\lambda z^3$. The Prym variety $P$ is
isogenous to a product of six elliptic curves:
$$
P \sim {\rm Jac}(C_{f_1})^2 \times {\rm Jac}(C_{f_3})^2 \times P^{\prime},
$$
where $P^{\prime}$ is a supersingular abelian surface whose trace of Frobenius
$t(P^{\prime})$ over $\F_{q}$  satisfies
$$
t(P^{\prime})=\cases{0 & if $m$ odd, \cr
- 2(q-1)+2t(C_{f_3})   & if $m$ even.\cr}
$$
\par
\noindent
{\sl Proof.} The splitting of ${\rm Jac}(C/\tau)$ follows directly from the
description of $C/\tau$ as a fibre product 
and the splitting ${\rm Jac}(C_{f_1+f_2})  \sim {\rm Jac}(C_{g_{\lambda}}) 
\times E$
as obtained in (2.4). Furthermore, using Theorem (2.3) we see that
$$
P \sim {\rm Jac}(C_{f_1}) \times {\rm Jac}(C_{f_2})\times {\rm Jac}(C_{f_3})
\times {\rm Jac}(C_{f_1+f_2+f_3})
$$
We know ${\rm Jac}(C_{f_1}) \cong {\rm Jac}(C_{f_2})$ and that 
${\rm Jac}(C_{f_1+f_2+f_3})$  splits up to isogeny as ${\rm Jac}(C_{f_3})$
and a $2$-dimensional factor $P^{\prime}$ 
which is supersingular and up to isogeny
a product of two elliptic curves. Using the map $x\mapsto w=x^3$ 
we see that $\# C_{f_1+f_2+f_3}(\F_{q})= \# C_{f_3}(\F_{q})$
if $m$ is odd which implies $t(P^{\prime})=0$, while for $m$ even
$$
\# C_{f_1+f_2+f_3}(\F_{q})-2=3( \# C_{f_3}(\F_{q})-2).
$$
This implies
$$
t(C_{f_1+f_2+f_3})-t(C_{f_3})=- 2(q-1) +2t(C_{f_3})
$$
and hence $t(P^{\prime})=- 2(q-1) +2t(C_{f_3})$. This proves the assertion.

\bigskip
\centerline{\bf\S 3. Bounds for $N(A,B)$}
\bigskip
\noindent
Since the curve $C=C_{A,B}$ has genus $13$ if $\lambda =A^2+A+1+B \neq 0$
the Hasse-Weil-Serre bound for the number of $\F_q$-rational
points $\# C_{A,B}(\F_q)$ says
$$
q+1-13[2\sqrt{q}] \leq \# C_{A,B}(\F_q) \leq q+1+13[2\sqrt{q}]. \eqno(12)
$$
The number $N(A,B)$ of $S_4$-orbits of solutions of (1)
with  distinct $x_i \in \F_q$ satisfies
$$
N(A,B)=(\# C_{A,B}(\F_q)-\hbox{\rm contribution of $x=0,1,\infty$})/24.
$$
It ${\rm Tr}(A+1)=0$ we have $12$ rational points in the fibres above 
$0,1,\infty$, while there are none if ${\rm Tr}(A+1)=1$. Then (12) implies
for $N(A,B)$ the inequalities
$$
(q-11-13[2\sqrt{q}])/24 \leq N(A,B) \leq (q+1+13[2\sqrt{q}])/24.
$$
By employing the decomposition of the Jacobian, especially  Corollary (2.4),
and taking
into account that the possible values of the trace of Frobenius $t$ of
supersingular elliptic curves  
are $t=0, \, \pm \sqrt{2q}$ for $q=2^m$ with $m$ odd we can refine
these bounds and 
we obtain our main result on the numbers $N(A,B)$: 

\proclaim (3.1) Theorem. For $q=2^m$ with $m$ odd the number 
$N(A,B)$ satisfies the following inequalities:
$$
(q-11-2\sqrt{2q}-8[2\sqrt{q}])/24 \leq N(A,B) \leq 
(q+1+2\sqrt{2q}+8[2\sqrt{q}] )/24. \eqno(13)
$$
\par

In the following table we illustrate this by listing the intervals
in which the numbers lie according to  (13). 
\smallskip
\noindent
{\bf Table 1}.
\smallskip
\vbox{
\bigskip\centerline{\def\quad{\hskip 0.6em\relax}
\def\quod{\hskip 0.5em\relax }
\vbox{\offinterlineskip
\hrule
\halign{&\vrule#&\strut\quod\hfil#\quad\cr
height2pt&\omit&&\omit&&\omit&&\omit&&\omit&&\omit&\cr
&$ q$&&$32$&&$128$&&$512$&&$2048$&&$8192$&\cr
height2pt&\omit&&\omit&&\omit&&\omit&&\omit&&\omit&\cr
\noalign{\hrule}
height2pt&\omit&&\omit&&\omit&&\omit&&\omit&&\omit&\cr
&${ \rm interval } $&&$[0,4]$&&$[0,14]$&&$[4,38]$&&$[50,120]$&&$[270,412]$&\cr
height2pt&\omit&&\omit&&\omit&&\omit&&\omit&&\omit&\cr
} \hrule}
}}
\bigskip

For some further reflections on $N(A,B)$ we restrict to the case
$q=2^m$ with $m$ odd. 
The practice of searching
for curves with many points tells us that is is highly improbable that
in a fibre product of curves the traces of Frobenius of the individual
components simultaneously reach their maximal (or minimal) value.
Hence it is very unlikely that the bounds given in (13) will be reached.

We intend to design an interval which contains
almost all values of $N(A,B)$ using the 
description of $C$ as a fibre product of the curves $C_{f_i}$ for $i=1,2,3$
given in (11) and 
a probabilistic  argument on the distribution of traces of Frobenius.

The curves $C_{f_1}$ and $C_{f_2}$ are supersingular elliptic curves
with the same trace of Frobenius $t=t(C_{f_i})=0,\, \pm \sqrt{2q}$.
The curve $C_{f_1+f_2}$ has genus $3$ 
which implies
$$
-3[2\sqrt{q}]\leq t(C_{f_1+f_2})\leq 3[2\sqrt{q}].
$$
So the trace of Frobenius for the normalization of the fibre product 
$C_{f_1}\times_{\P^1 } C_{f_2}$ satisfies
$$
-3[2\sqrt{q}]-2\sqrt{2q} \leq t \leq 3[2\sqrt{q}] +2\sqrt{2q}.
$$
We compute bounds for the number of $x \in \P^1-\{ 0, 1, \infty\}$
above which we have $4$ points in the fibre of $C_{f_1}\times_{\P^1} C_{f_2}$
If ${\rm Tr}(A+1)=0$ we find in total $8$ points
above $x=0,1,\infty$, while we find none if ${\rm Tr}(A+1)=1$.
Subsequently we take into account that completely splitting 
$x\in \P^1-\{0,1,\infty\}$ occur in pairs $(x,1/x)$ 
and we obtain the following Proposition.

\proclaim (3.2) Proposition. If we let
$$
M(f_1,f_2)= {1 \over 2}  \# 
\{ x \in \P^1(\F_q) -\{0,1, \infty\} 
\hbox{ \rm :  $x$ splits completely in $C_{f_1}\times_{\P^1} C_{f_2}$} \}
$$
then we have for ${\rm Tr}(A+1)=0$
$$
{q-7-3[2\sqrt{q}]-2\sqrt{2q} \over 8}  \leq
M(f_1,f_2) \leq
{q-7+3[2\sqrt{q}]+2\sqrt{2q} \over 8} 
$$
and for ${\rm Tr}(A+1)=1$
$$
{q+1-3[2\sqrt{q}]-2\sqrt{2q} \over 8}  \leq
M(f_1,f_2) \leq
{q+1+3[2\sqrt{q}]+2\sqrt{2q} \over 8}.
$$
\par
\bigskip
We now consider the effect of the elliptic curve $C_{f_3}$ in the fibre product. The $j$-invariant of $C_{f_3}$ is $\lambda^{-4} \in \F_q^*$. This implies that
$t(C_{f_3})$ is odd. For ${\rm Tr}(A+1)=0$ we have 
$t(C_{f_3})\equiv 1 (\bmod \, 4)$ and there are $4$ 
rational points together above $x=0,1,\infty$, while if ${\rm Tr}(A+1)=1$ 
we have $t(C_{f_3})\equiv 3 (\bmod \, 4)$ and $2$ rational points 
above $0,1,\infty$. Furthermore, each element of $\F_q^*$ occurs exactly once
as $j$-invariant in the family of curves $C_{f_3}$. That implies that
$t(C_{f_3})$ assumes each odd integer value in the interval
$[-[2\sqrt{q}],[2\sqrt{q}]]$. So the number of completely splitting pairs
assumes each integral value in the intervals mentioned in the
following Proposition.

\proclaim (3.3) Proposition. If we let
$$
M(f_3)= {1 \over 2} 
\# \{x \in \P^1(\F_q)-\{0,1,\infty\} :
\hbox{\rm $x$  splits completely in $C_{f_3}$}\}$$
then $M(f_3)$ assumes all integer
values in the intervals
$$
\eqalign{
[{q-3-[2\sqrt{q}]\over 4}, {q-3+[2\sqrt{q}]\over 4}] & \quad {\rm if} \quad {\rm Tr}(A+1)=0, \cr
[{q-1-[2\sqrt{q}]\over 4}, {q-1+[2\sqrt{q}]\over 4}] & \quad {\rm if} \quad {\rm
 Tr}(A+1)=1. \cr}
$$
\par
\medskip
Finally, we combine the two preceding propositions via a heuristic argument.
Let 
$$
\eqalign{M(f_1,f_2,f_3)= &\cr
{1 \over 2}  \# \{ x \in & \P^1(\F_q)-  \{0,1,\infty\}: \hbox{
\rm $x$ splits completely on 
$C_{f_1}\times_{\P^1} C_{f_2}\times_{\P^1} C_{f_3}$}\}\cr}
$$
Since there are $(q-2)/2$ pairs $(x,1/x)$ ($x \neq 0,1,\infty$) we expect
$$
\eqalign{
2\big( {q-3-[2\sqrt{q}] \over 4}\big )\big( {q-7-3[2\sqrt{q}]
-2\sqrt{2q} \over 8} \big) / (q-2) & \leq M(f_1,f_2,f_3) \leq \cr
2\big( {q-1+[2\sqrt{q}] \over 4}\big) & \big(
{q+1+3[2\sqrt{q}]+2\sqrt{2q} \over 8
} \big)  / (q-2). \cr }
$$
If we work this out and neglect terms of order $1/\sqrt{q}$ and lower we
find
$$
{ q-4[2\sqrt{q}] -2\sqrt{2q} +4+4\sqrt{2} \over 16}
\leq M(f_1,f_2,f_3) \leq
{q+4[2\sqrt{q}] +2\sqrt{2q} +14 +4\sqrt{2} \over 16} \eqno(14)
$$
Each completely splitting pair yields $16$ solutions of (1) 
so to estimate the number of $S_4$-orbits of solutions $N(A,B)$ 
we multiply the interval by $16/24$ to get an interval $I$.
Since $N(A,B)$ is even we adapt the endpoints of the interval 
$I$ just obtained slightly.
Namely we consider the smallest interval with endpoints the positive even 
integers which contains $I$ and we 
denote this interval by $I^{\rm even}$. 

\proclaim (3.4) Heuristics. 
The odds are that the values of $N(A,B)$ are in the interval
$$
I^{\rm even}=[{q-4[2\sqrt{q}]-2\sqrt{2q}+4\sqrt{2}+4 \over 24},
{q+4[2\sqrt{q}]+2\sqrt{2q}+4\sqrt{2} +14 \over 24}]^{\rm even}.
$$
\par

\noindent
We illustrate this by a little table.
\smallskip
\noindent
{\bf Table 2.}
\smallskip
\vbox{
\bigskip\centerline{\def\quad{\hskip 0.6em\relax}
\def\quod{\hskip 0.5em\relax }
\vbox{\offinterlineskip
\hrule
\halign{&\vrule#&\strut\quod\hfil#\quad\cr
height2pt&\omit&&\omit&&\omit&&\omit&&\omit&&\omit&\cr
&$q$&&$32$&&$128$&&$512$&&$2048$&&$8192$&\cr
height2pt&\omit&&\omit&&\omit&&\omit&&\omit&&\omit&\cr
\noalign{\hrule}
height2pt&\omit&&\omit&&\omit&&\omit&&\omit&&\omit&\cr
&${ \rm interval } $&&$[0,6]$&&$[0,12]$&&$[10,34]$&&$[64,108]$&&$[300,384]$&\cr
height2pt&\omit&&\omit&&\omit&&\omit&&\omit&&\omit&\cr
} \hrule}
}}
\bigskip

\bigskip
\centerline{\bf \S 4. Numerical results}
\bigskip
\noindent
In order to obtain numerical results on $N(A,B)$ to test our heuristics
the first remark is that $N(A_1,B)=N(A_2,B)$ if ${\rm Tr}(A_1)={\rm Tr}(A_2)$.
So we have to distinguish only between ${\rm Tr}(A)=0$ and ${\rm Tr}(A)=1$.
We shall compute the trace of Frobenius for the seven factors of our
Jacobian.  We shall write $f_4=f_1+f_2$, $f_5=f_1+f_3$, $f_6=f_2+f_3$ and 
$f_7=f_1+f_2+f_3$. The Jacobians of the 
curves $C_{f_i}$ given by $y^2+y=f_i$ for $i=1,\ldots,7$
constitute the seven factors of ${\rm Jac}(C_{A,B})$. We write
$$
n_{f_i}= \# \{ x \in \F_q^* : {\rm Tr}(f_i(x))=0 \}.
$$

\proclaim (4.1)
Proposition. The number of solutions $N(A,B)$ over $\F_{q=2^m}$
with $m$ odd of the system (1) with $\lambda = A^2+A+1+B\neq 0$ is given by
$$
N(A,B)= { 2q-2 -2\{ n_{f_1}+n_{f_2}+n_{f_3}-n_{f_4}-n_{f_5}-n_{f_6}+n_{f_7} \}
\over 24} \qquad \hbox{\rm if ${\rm Tr}(A)=0$},
$$
and
$$
N(A,B)={ -6q-2+2\sum_{i=1}^7 n_{f_i} \over 24} \qquad \hbox{\rm if
${\rm Tr}(A) =1$}.
$$
\par
\noindent
{\sl Proof.} As just explained we may take $A=0$ or $A=1$. 
Then $\lambda= B+1\neq 0$
and we set $f_1=(B+1)(x^3+x)$, $f_2=(B+1)(1/x^3+1/x)$ and $f_3=(B+1)(x+1/x)$.
Then $C_{1,B}=C_{f_1}\times_{\P^1} C_{f_2}\times_{\P^1} C_{f_3}$
and $C_{0,B}=C_{f_1+1}\times_{\P^1} C_{f_2+1}\times_{\P^1} C_{f_3+1}$.
As in Theorem (2.3) the curves $C_{f_i}$ for $i=4,\ldots,7$ give the 
remaining traces of Frobenius.

The trace of Frobenius $t(C_{f_i})$ is of the form
$$
t(C_{f_i})= q+1-2n_{f_i} - a_i,
$$
where $a_i$ is the contribution of $x=0, \infty$, while the trace of
Frobenius of $C_{f_i+1}$ is
$$
t(C_{f_i})= -q+3+2n_{f_i} - b_i ,
$$
where $b_i$ is the contribution of $x=0, \infty$. By analyzing these
contributions from $0$ and $\infty$ one gets the Proposition.
\bigskip
We now give tables with the distribution of the numbers $N(A,B)$ for
$q=2^m$ with $m$ odd and $5 \leq m \leq 13$. These tables are obtained
by computing the numbers $n_{f_i}$ and they solve the coset weight
distribution problem for the corresponding $BCH(3)$ codes. The first
unknown case up to now was $q=2^9$, see [C-Z]. Moreover, the tables
confirm our heuristics. We list the frequencies divided by $q/2$.
\medskip\noindent
{\bf Table 3}.
\medskip
$q=2^5$
\smallskip
\vbox{
\bigskip\centerline{\def\quad{\hskip 0.6em\relax}
\def\quod{\hskip 0.5em\relax }
\vbox{\offinterlineskip
\hrule
\halign{&\vrule#&\strut\quod\hfil#\quad\cr
height2pt&\omit&&\omit&&\omit&\cr
&$N(A,B)$&&$0$&&$2$&\cr
height2pt&\omit&&\omit&&\omit&\cr
\noalign{\hrule}
height2pt&\omit&&\omit&&\omit&\cr
&${\rm frequency}$&&$27$&&$35$&\cr
height2pt&\omit&&\omit&&\omit&\cr
} \hrule}
}}
\bigskip

$q=2^7$

\smallskip
\vbox{
\bigskip\centerline{\def\quad{\hskip 0.6em\relax}
\def\quod{\hskip 0.5em\relax }
\vbox{\offinterlineskip
\hrule
\halign{&\vrule#&\strut\quod\hfil#\quad\cr
height2pt&\omit&&\omit&&\omit&&\omit&&\omit&&\omit&&\omit&\cr
&$N(A,B)$&&$0$&&$2$&&$4$&&$6$&&$8$&&$10$&\cr
height2pt&\omit&&\omit&&\omit&&\omit&&\omit&&\omit&&\omit&\cr
\noalign{\hrule}
height2pt&\omit&&\omit&&\omit&&\omit&&\omit&&\omit&&\omit&\cr
&${ \rm frequency} $&&$2$&&$28$&&$98$&&$84$&&$35$&&$7$&\cr
height2pt&\omit&&\omit&&\omit&&\omit&&\omit&&\omit&&\omit&\cr
} \hrule}
}}
\bigskip

$q=2^9$

\smallskip
\vbox{
\bigskip\centerline{\def\quad{\hskip 0.6em\relax}
\def\quod{\hskip 0.5em\relax }
\vbox{\offinterlineskip
\hrule
\halign{&\vrule#&\strut\quod\hfil#\quad\cr
height2pt&\omit&&\omit&&\omit&&\omit&&\omit&&\omit&&\omit&&\omit&&\omit&&\omit&&\omit&&\omit&\cr
&$N(A,B)$&&$12$&&$14$&&$16$&&$18$&&$20$&&$22$&&$24$&&$26$&&$28$&&$30$&&$32$&\cr
height2pt&\omit&&\omit&&\omit&&\omit&&\omit&&\omit&&\omit&&\omit&&\omit&&\omit&&\omit&&\omit&\cr
\noalign{\hrule}
height2pt&\omit&&\omit&&\omit&&\omit&&\omit&&\omit&&\omit&&\omit&&\omit&&\omit&&\omit&&\omit&\cr
&${ \rm frequency} $&&$18$&&$21$&&$117$&&$180$&&$148$&&$195$&&$199$&&$81$&&$36$&&$18$&&$9$&\cr
height2pt&\omit&&\omit&&\omit&&\omit&&\omit&&\omit&&\omit&&\omit&&\omit&&\omit&&\omit&&\omit&\cr
} \hrule}
}}
\bigskip

$q=2^{11}$

\smallskip
\vbox{
\bigskip\centerline{\def\quad{\hskip 0.6em\relax}
\def\quod{\hskip 0.5em\relax }
\vbox{\offinterlineskip
\hrule
\halign{&\vrule#&\strut\quod\hfil#\quad\cr
height2pt&\omit&&\omit&&\omit&&\omit&&\omit&&\omit&&\omit&&\omit&&\omit&&\omit&&
\omit&&\omit&\cr
&$N(A,B)$&&$66$&&$68$&&$70$&&$72$&&$74$&&$76$&&$78$&&$80$&&$82$&&$84$&&$86$&\cr
height2pt&\omit&&\omit&&\omit&&\omit&&\omit&&\omit&&\omit&&\omit&&\omit&&\omit&&
\omit&&\omit&\cr
\noalign{\hrule}
height2pt&\omit&&\omit&&\omit&&\omit&&\omit&&\omit&&\omit&&\omit&&\omit&&\omit&&
\omit&&\omit&\cr
&${ \rm frequency} $&&$22$&&$66$&&$88$&&$55$&&$176$&&$264$&&$187$&&$374$&&$374$&
&$374$&&$451$&\cr
height2pt&\omit&&\omit&&\omit&&\omit&&\omit&&\omit&&\omit&&\omit&&\omit&&\omit&&
\omit&&\omit&\cr
\noalign{\hrule}
height2pt&\omit&&\omit&&\omit&&\omit&&\omit&&\omit&&\omit&&\omit&&\omit&&\omit&&
\omit&&\omit&\cr
&$N(A,B)$&&$88$&&$90$&&$92$&&$94$&&$96$&&$98$&&$100$&&$102$&&$104$&&$106$&&108&\cr
height2pt&\omit&&\omit&&\omit&&\omit&&\omit&&\omit&&\omit&&\omit&&\omit&&\omit&&
\omit&&\omit&\cr
\noalign{\hrule}
height2pt&\omit&&\omit&&\omit&&\omit&&\omit&&\omit&&\omit&&\omit&&\omit&&\omit&&
\omit&&\omit&\cr
&${ \rm frequency} $&&$365$&&$341$&&$275$&&$341$&&$154$&&$44$&&$55$&&$33$&&$11$&
&$22$&&$22$&\cr
height2pt&\omit&&\omit&&\omit&&\omit&&\omit&&\omit&&\omit&&\omit&&\omit&&\omit&&
\omit&&\omit&\cr
} \hrule}
}}
\bigskip
$q=2^{13}$
\smallskip
In this case we encounter a new phenomenon.
The function $N(A,B)$ assumes even values in the interval $[290,390]$,
but not all even values are taken.  This  contradicts
the expectation of [C-Z] that the values form a sequence of
even integers without gaps.
The frequency divided by $q/2$ 
of the value $290+2\ell$ with $0 \leq \ell \leq 50$ is given by
$$
13 \, \gamma_{\ell} + \cases{1 & if $\ell=11$,\cr
1 & if $\ell= 37$, \cr
0 & else, \cr } 
$$
where $\gamma =(\gamma_0,\ldots,\gamma_{50})$ is the vector
$$
\eqalign{
\gamma =(1,0,1,&0,1,0,6,3,5,5,12,7,19,15,22,25,37,40,43,37,35,60,
54,72,72,58,65,\cr 
&61,57,57,63,48,35,44,34,34,25,29,25,15,9,7,2,3,7,
3,3,1,0,1,2). \cr}
$$
In accordance with our heuristics
less than $1 \, \%$ of the $N(A,B)$ lie outside the interval
$[300,384]$.

\bigskip
\centerline{\bf \S5. The covering radius}
\bigskip
\noindent
A problem in coding theory that precedes the coset weight distribution
problem is the determination of the covering radius. It is defined for
a binary linear code $C$ of length $n$ as the smallest integer $\rho$
such that the spheres of radius $\rho$ around the codewords cover $\F_2^n$.
Equivalently, it is the maximum weight of a coset leader (by which we mean a
vector of minimum weight in a coset of $C$ in $\F_2^n$). It is an interesting
parameter of a code since it provides information on the performance
of the code when used in data compression.

In a series of papers [H-B],[A-M] and [H], of which [H-B] and [H] 
treat the case  $m$ even and [A-M] the case $m$ odd,
  it was proved that the $BCH(3)$ code
of length $n=2^m-1$ has covering radius
$$
\rho(BCH(3)) =5 \qquad {\rm for} \quad m\geq 4.
$$
The proofs for the various cases are very different. Using algebraic geometry
we can give a unified proof.

In order to prove that $\rho (BCH (3)) = 5$ we have to show that for 
every $(A,B,C) \in \F_q^3$ the system of equations:

$$
\eqalign{
&x_1 + \ldots + x_5 = A, \cr
&x_{1}^3 + \ldots +x_{5}^3 = B, \cr
&x_{1}^{5} + \ldots +x_{5}^{5} = C,\cr}\eqno(15)
$$
has a solution $(x_1,..., x_5) \in \F_{q}^5$ for any 
$(A, B, C) \in \F _{q}^3$. By replacing $x_i$ by $x_i+A$ 
we may assume without loss of generality that $A = 0$ and 
$(B, C) \not= (0, 0)$. If we then homogenize (15) the system
$$
\sum_{i = 1}^{5} x_i = 0, ~~\sum_{i = 1}^{5} x_{i}^3 = Bx_{0}^3,~~\sum_{i = 1}
^{5} x_{i}^5 = C x _{0}^5 \eqno (16)
$$
defines a projective variety $V$ of dimension $2$ in the five dimensional
projective space $\P^5$.

We intersect $V$ with the hyperplane $x_0+x_5=0$ and obtain a system
of equations of the form (2). By using the results of Section 1 (especially
Corollary  (1.3)) one can easily show that $\rho(BCH(3)) =5 $ for $m\geq 10$.
We leave the details to the reader.
\bigskip
As a final remark we would like to point out that we think that
many more problems on cyclic codes can be attacked succesfully 
using methods from algebraic geometry as is done in this paper. 
We refer to [C] for a list of such problems.
\bigskip
\centerline{\bf References}
\bigskip
\noindent
[A-M] E.F.\ Assmus, Jr., H.F.\ Mattson, Jr.: Some $3$-error-correcting BCH
codes have covering radius $5$. 
{\sl IEEE Trans.\ Info.\ Th., \bf 22}, 1976, p.\ 348--349.
\smallskip
\noindent
[C] P.\ Charpin: Open problems on cyclic codes. Handbook of Coding
Theory I. V.S.\ Pless, W.C.\ Huffman Eds., Elsevier Science BV,
Amsterdam, 1998, p.\ 963--1063.
\smallskip
\noindent
[C-Z] P.\ Charpin, V.\ Zinoviev: 
On coset weight distributions of the $3$-error-correcting BCH codes. 
{\sl SIAM J. Discrete Math.\ \bf  10}, (1997),  128--145. 
\smallskip
\noindent
[H]  T.\ Helleseth, All binary 3-error-correcting BCH codes of length
$2^m-1$  have covering radius $5$. 
{\sl IEEE Trans.\ Info.\ Th.\ \bf 24}, 1978, p.\  257--258.
\smallskip
\noindent
[H-B] J.A.\ van der Horst, T.\ Berger: Complete decoding of 
triple-error-correcting binary BCH codes. 
{\sl IEEE Trans.\ Info.\ Th.\  \bf 22}, 1976, p.\ 138--147.
\smallskip
\noindent
\bigskip
\bigskip
\settabs3 \columns
\+G. van der Geer  &&M. van der Vlugt\cr
\+Faculteit
Wiskunde en Informatica &&Mathematisch Instituut\cr
\+Universiteit van
Amsterdam &&Rijksuniversiteit te Leiden \cr
\+Plantage Muidergracht 24&&Niels Bohrweg 1 \cr
\+1018 TV Amsterdam
&&2333 CA Leiden \cr
\+The Netherlands &&The Netherlands \cr
\+{\tt geer@science.uva.nl} &&{\tt vlugt@math.leidenuniv.nl} \cr

\end